\documentclass[sort&compress,3p]{elsarticle}

\usepackage{amsfonts}
\usepackage{amssymb}
\usepackage{amsthm}
\usepackage{amsmath}
\usepackage{hyperref}
\usepackage{enumitem}

\newtheorem{theorem}{Theorem}

\journal{Applied Mathematics Letters}

\begin{document}

\begin{frontmatter}

\title{Stability of Equilibria in Modified-Gradient Systems}

\author[bjr]{Benjamin J. Ridenhour\corref{cor1}}
\ead{bridenhour@uidaho.edu}
\author[jrr]{Jerry R. Ridenhour}

\address[bjr]{Department of Biological Sciences, University of Idaho, Moscow, ID 83844, USA}
\address[jrr]{Department of Mathematics and Statistics, Utah State University, Logan, UT 84322, USA}

\cortext[cor1]{Corresponding author. Phone: +1 208 885 7231. Fax: +1 208 885 7905}

\begin{abstract}
Motivated by questions in biology, we investigate the stability of equilibria
of the dynamical system $\mathbf{x}^{\prime}=P(t)\nabla f(x)$ which arise as
critical points of $f$, under the assumption that $P(t)$ is positive
semi-definite. It is shown that the condition $\int^{\infty}\lambda
_{1}(P(t))~dt=\infty$, where $\lambda_{1}(P(t))$ is the smallest eigenvalue of
$P(t)$, plays a key role in guaranteeing uniform asymptotic stability and in
providing information on the basis of attraction of those equilibria.
\end{abstract}


\end{frontmatter}

%\linenumbers

\section{Introduction}
The evolution of continuous phenotypes, for example height, by means of natural
selection is a central theme in evolutionary biology. The breeder's equation
$(R = h^2 s)$ was first introduced by Lush in 1937 \cite{Lush1937} to predict the
change in phenotype $(R)$ with respect to the heritability $(h^2)$ and strength of
natural selection $(s)$. In a seminal series of papers, the breeder's equation was updated
to the so-called multivariate breeder's equation by Lande 
\cite{lande_natural-selection_1976,lande_quantitative_1979} and Lande
and Arnold \cite{lande_measurement_1983}. The multivariate breeder's equation
is often presented in varying forms such as $\Delta \bar{z}(t) = h^2 \sigma^2 \partial 
\mbox{\,ln}(\bar{W})/\partial \bar{z}(t)$ \cite{lande_natural-selection_1976}, 
$\Delta \bar{\mathbf{z}} = \mathbf{G} \nabla \mbox{ln}(\bar{\mathbf{W}})$
\cite{lande_quantitative_1979}, $\Delta \bar{\mathbf{z}} = \mathbf{G P^{-1} s}$ 
\cite{lande_measurement_1983}, and $\Delta \bar{\mathbf{z}} = \mathbf{G \beta}$ 
\cite{Blows:2014jx}, as well as continuous-time counterparts (i.e., $\mbox{d}\mathbf{\bar{z}}/\mbox{d}t$);
all of these forms reduce to the concept that the change in mean phenotype $(\bar{\mathbf{z}})$
over time is given by the product of a genetic variance-covariance matrix $(\mathbf{G})$
and the gradient of the logarithm of the average fitness function $(\bar{W}(\bar{\mathbf{z}}))$.
As of December 2015, Web of Science indicates that the papers by Lande 
\cite{lande_natural-selection_1976,lande_quantitative_1979} and Lande
and Arnold \cite{lande_measurement_1983} have garnered at least 791, 1\,442,
and 2\,852 citations, respectively, which gives some idea of the impact these 
works have had on evolutionary biology and related fields.

One of the critical assumptions in much of this research is that the so-called 
$\mathbf{G}$-matrix is constant. A Web of Science search indicates at least 
175 papers on the constancy and form of the $\mathbf{G}$-matrix with 66 of those 
published since January 2010 (a broader search on ``genetic constraints'' reveals 
many more relevant publications). The principal concern is that the  $\mathbf{G}$-matrix 
limits how traits evolve and approach their evolutionary optima 
\cite{Arnold:2008bf,Blows:2005um,lande_quantitative_1979}. For example, Dickerson 
\cite{Dickerson55} studied a special case of equal genetic variances which produces a 
$\mathbf{G}$-matrix with a zero eigenvalue, thus preventing evolution along some
trajectories. Furthermore, Pease and Bull \cite{Pease:1988ea} examined ``ill-conditioned'' 
$\mathbf{G}$-matrices where the ratio of the largest to the smallest eigenvalue is large and concluded that 
the speed of evolution toward an optimum is greatly reduced.
Other work has suggested that the number of dimensions in the system affects stability \cite{Debarre:2014kf}.
However, formal criteria for when and how an evolutionary system will converge upon an 
equilibrium are lacking. While most research considers
$\mathbf{G}$ to be constant, it is widely recognized that $\mathbf{G}$
itself is expected to evolve over time \cite{Arnold:2008bf,Bjorklund:2015hn,
lande_measurement_1983, Laughlin:2015kj}. Considering $\mathbf{G}$  to
be time-varying further muddies the waters of whether such systems
approach and are stable at existing equilibria and lacks formal mathematical
treatment. 

Similarly, considerable interest has been paid to rugged fitness landscapes 
where the average fitness function has multiple peaks
(optima) \cite{Burger:1986kv,Nahum:2015dz,Neidhart:2014dr,Whitlock:1995ec}.
Exploration of fitness landscapes, in other words movement between different
optima, is a key part of Wright's shifting balance theory \cite{wright_evolution_1931}.
Despite interest in which evolutionary optimum the population mean phenotype
will evolve toward, conventional wisdom that the nearest optimum is favored
or numerical methods are relied upon. In fact, some research has shown that the nearest 
optimum is not always the one favored by evolution \cite{Burger:1986kv}. As in 
the case of stability analyses, no rigorous analysis of if and when a particular
optimum will be evolved toward has been performed.

The contributions of this paper are threefold. First, we rigorously analyze
the modified gradient system commonly used to model the evolution of continuous traits
for the existence and stability of equilibria. Our research shows that biologists can simply search for the isolated local maxima of a fitness function; these points are guaranteed to be
at least uniformly stable. Second, in cases where the smallest eigenvalue
of the G-matrix, $\lambda_{1}(P(t))$ in our notation, meets the condition $\int_{0}^{\infty
}\lambda_{1}(P(t))dt=\infty$, then the equilibrium is guaranteed to be
uniformly asymptotically stable. Finally, an understanding of the inverse
image under the fitness function $f$ of intervals of the form $(c,f(\overline
{x}))$ gives information on the basin of attraction of an equilibrium at
$\overline{x}$. Taken together, these contributions greatly enhance our ability to analyze and understand multivariate phenotypic evolution.

\section{Existence and Stability of Equilibria}

Let $x=(x_{1},x_{2},\ldots,x_{n})$ denote a point in $\mathbb{R}^{n}$, and let
$\mathbf{x}=[x_{1}\hspace*{0.05in}x_{2}\hspace*{0.05in}\cdots\hspace
*{0.05in}x_{n}]^{T}$ be the corresponding $n\times1$ vector equivalent. We use
the Euclidean norm as a measure of distance and we let $B_{\delta}%
(\overline{x})$ denote the open ball of radius $\delta$ centered at
$\overline{x}$. The object is to determine the stability of equilibrium
solutions of the $n$-dimensional modified-gradient system%

\begin{equation}
\mathbf{x}^{\prime}=P(t)\mathbf{\nabla}f(x). \tag{E}%
\end{equation}

\noindent Note that the continuous-time multivariate breeder's equation is of this form with
$P(t)$ being the time-dependent $\mathbf{G}$-matrix, and $f$ being 
$\mbox{ln}\,\bar{W}$. We assume throughout that the following 
hypothesis holds:

\textit{H}$_{0}$\textit{: }$D$\textit{ is a domain in }$R^{n}$\textit{, }%
$f$\textit{ is a real-valued }$C_{1}$\textit{ (i.e., continuous with
continuous partials) function defined on }$D$\textit{, }$t$\textit{ is
nonnegative, the gradient of }$f$\textit{ denoted by }$\nabla f$\textit{ has
components which are }$C_{1}$\textit{ on }$D,$\textit{ and }$P(t)$\textit{ is
an }$n\times n$\textit{ matrix-valued function with }$C_{1}$\textit{-entries
that is defined and positive semi-definite for }$t\geq0$\textit{.}

\textit{H}$_{0}$ guarantees that, for any $t_{0}\geq0$ and any $x_{0}$ in $D$,
there is a unique solution of (E) satisfying the initial condition
$x(t_{0})=x_{0}$. The assumption that $P(t)$ is positive semi-definite is 
consistent with biological applications because the $\mathbf{G}$-matrix is a 
variance-covariance matrix, and variance-covariance matrices are always symmetric,
positive semi-definite matrices. 

If $f$ has an isolated maximum value at a point $x=\overline{x}$ of $D$, then
we know from calculus that $\nabla f(\overline{x})=\mathbf{0}$ so
$\mathbf{x}=\overline{\mathbf{x}}$ is an equilibrium (i.e., constant in time)
solution of (E). We investigate the stability of such equilibria. Although a
translation always allows one to assume the equilibrium point is at
$\mathbf{x}=\mathbf{0}$, we will continue to assume, because of our interest
in evolutionary applications, that $\mathbf{x}=\overline{\mathbf{x}}$ is the
equilibrium solution.

We let $\lambda_{1}(P(t))$ denote the smallest eigenvalue of $P(t)$ and
introduce the eigenvalue condition
\begin{equation}
\int_{0}^{\infty}\lambda_{1}(P(t))dt=\infty. \tag{EC}%
\end{equation}
This condition will play an important role in what follows.

We follow the definitions of uniform stability and uniform asymptotic
stability as given by Hartman \cite{hartman_ordinary_1964}. In contrast to some definitions,
this definition of uniform asymptotic stability gives uniformity in the choice
of starting time $t_{0}$ and does not involve the rate at which solutions tend
to the equilibrium solution. Consider the following additional hypotheses:

\textit{H}$_{1}$\textit{:~ }$f$\textit{ has an isolated local maximum value at
the point }$\overline{x}\in$\textit{ }$D$\textit{;}

\textit{H}$_{2}$\textit{:~ }$\overline{x}$\textit{ is an isolated critical
point of }$f;$\textit{ and}

\textit{H}$_{3}$\textit{:~ eigenvalue condition (EC) holds.}

Our stability results are contained in the following theorem. 

\begin{theorem}[Stability and Asymptotic Stability]\label{Th 1} $\,$ 

 	 \begin{enumerate}[label={(\roman*)}, itemindent = 1em]
		\item If H$_{0}$ and
H$_{1}$ hold, then $x=\overline{x}$ is a uniformly stable equilibrium solution
of $(E)$.
		\item If H$_{0}$, H$_{1}$, H$_{2}$ and H$_{3}$ all hold, then
$x=\overline{x}$ is a uniformly asymptotically stable equilibrium solution of
$(E)$.
	\end{enumerate}

\begin{proof}
Suppose H$_{0}$ and H$_{1}$ hold. Let $M=f(\overline{x})$ and define the
function $V$ on $D$ by $V(x)=M-f(x)$. For a solution $x(t)$ of $(E)$, let
$V_{x}$ be the function defined by $V_{x}(t)=V(x(t))$ for $t$ in the interval
of existence of $x(t)$. The so-called trajectory derivative is then given by
\begin{align*}
V_{\mathbf{x}}^{\prime}(t)  &  =\dfrac{\partial V}{\partial x_{1}}%
(x(t))x_{1}^{\prime}(t)+\dfrac{\partial V}{\partial x_{2}}(x(t))x_{2}^{\prime
}(t)+\cdots+\dfrac{\partial V}{\partial x_{n}}(x(t))x_{n}^{\prime}(t)\\
&  =\mathbf{x}^{\prime}(t)\cdot\mathbf{\nabla}V(x(t))=-(P(t)\mathbf{\nabla
}f(x(t)))\cdot\mathbf{\nabla}f(x(t)).
\end{align*}
So, $V(\overline{x})=0$ and, in an appropriately chosen neighborhood of
$\overline{x}$, $V(x)>0$ for $x\neq\overline{x}$ and the trajectory
derivatives are nonpositive since $P(t)$ is positive semi-definite. By
standard Lyapunov theory (see, for example, Theorem 8.3, p. 40 of
Hartman \cite{hartman_ordinary_1964}), $x=\overline{x}$ is a uniformly stable equilibrium
solution of $(E)$. This proves (i).\newline\hspace*{0.25in}Now assume H$_{0}$,
H$_{1}$, H$_{2}$ and H$_{3}$ all hold. Let $M$ and $V$ be defined as in the
proof of (i). Suppose $\varepsilon>0$ is given. Since $\overline{x}$ is an
interior point of $D$ and since H$_{1}$ and H$_{2}$ hold, we can restrict
$\varepsilon>0$ to be so small that $x\in D$, $f(x)<M$ and $\mathbf{\nabla
}f(x)\neq\mathbf{0}$ for $0<\left\vert x-\overline{x}\right\vert
\leq\varepsilon$. Since we have uniform stability by Part (i), we find
$\delta>0$ with $\delta<\varepsilon$ such that any solution $x(t)$ satisfying
$\left\vert x(t_{0})-\overline{x}\right\vert <\delta$ at some time $t_{0}%
\geq0$ also satisfies $\left\vert x(t)-\overline{x}\right\vert <\varepsilon$
for all $t\geq t_{0}$.\newline\hspace*{0.25in}Let $x(t)$ be any solution with
$\left\vert x(t_{0})-\overline{x}\right\vert <\delta$ at some time $t_{0}%
\geq0$. To complete the proof of Part (ii), we need to show that
$\lim_{t\rightarrow\infty}x(t)=\overline{x}$. Since we have uniqueness of
solutions to initial value problems, we can assume $x(t)\neq\overline{x}$ for
$t\geq t_{0}$. Also, for $t\geq t_{0}$, $P(t)$ is positive semi-definite and
$\left\vert x(t)-\overline{x}\right\vert <\varepsilon$ so we have
$V_{\mathbf{x}}(t)=M-f(x(t))>0$ and $V_{\mathbf{x}}^{\prime}%
(t)=-(P(t)\mathbf{\nabla}f(x(t)))\cdot\mathbf{\nabla}f(x(t))\leq0$. Therefore,
$c=\lim_{t\rightarrow\infty}V_{\mathbf{x}}(t)$ exists with $c\geq0$. We prove
that $c=0$. Suppose not, then $c>0$. Since $V(\overline{x})=0$, we use the
continuity of $V$ to choose $\delta_{1}$ with $0<\delta_{1}<\delta$ so that
$V(x)<c$ for $|x-\overline{x}|<\delta_{1}.$ Because $V_{\mathbf{x}%
}(t)=V(x(t))>c$ for all $t\geq t_{0,}$ the trajectory $x(t)$ stays in the
region $\{x:\delta_{1}\leq\left\vert x-\overline{x}\right\vert \leq
\varepsilon\}$ for all $t\geq t_{0}.$ On this compact region, the continuous
function $\nabla f(x)\cdot\nabla f(x)$ is positive and hence assumes a
positive minimum value $m_{1}$ at some point in the set. Hence, for $t\geq
t_{0}$, we get that
\[
V_{x}^{^{\prime}}(t)\leq-\lambda_{1}(P(t))\nabla f(x(t))\cdot\nabla
f(x(t))\leq-m_{1}\lambda_{1}(P(t)).
\]
Consequently,
\[
V_{x}(t)=V_{x}(t_{0})+\int_{t_{0}}^{t}V_{x}^{^{\prime}}(s)ds\leq V_{x}%
(t_{0})-\int_{t_{0}}^{t}m_{1}\lambda_{1}(P(s))ds
\]
for $t\geq t_{0}$. But then (EC) implies that $V_{x}(t)\rightarrow-\infty$ as
$t\rightarrow\infty$ contradicting that $V_{x}(t)$ stays positve. This proves
that $\lim_{t\rightarrow\infty}V_{x}(t)=0$.\newline\hspace*{0.25in}Finally, we
prove that $\lim_{t\rightarrow\infty}x(t)=\overline{x}$. Suppose not. Then
there exists $\varepsilon_{1}$ with $0<\varepsilon_{1}<\varepsilon$ and
arbitrarily large values of $t$ where $\left\vert x(t)-\overline{x}\right\vert
\geq\varepsilon_{1}$. For such $t$, $V_{x}(t)=V\left(  x(t)\right)  \geq
m_{2}$ where $m_{2}$ is defined to be the minimum value of $V(x)$ on the
compact set $\{x:\varepsilon_{1}\leq x\leq\varepsilon\}$. Because $m_{2}>0$,
this contradicts $\lim_{t\rightarrow\infty}V_{x}(t)=0$ completing the proof.
\end{proof}
\end{theorem}

\emph{Example 2.1 Asymptotic Stability requires Eigenvalue Condition.}
The following example illustrates that an eigenvalue condition, such as we
have given in (EC), is necessary in order to obtain asymptotic stability.
Let $P(t)=\left[
\begin{array}
[c]{cc}%
(t+1)^{-2} & 0\\
0 & (t+1)^{-1}%
\end{array}
\right]  $ and $f(x_{1},x_{2})=4-(x_{1}-1)^{2}-(x_{2}-1)^{2}$ and consider the associated 
system $\mathbf{x}^{\prime}=P(t)\nabla f(x)$ for $t \geq 0$.
%\[
%\mathbf{x}^{\prime}=\left[
%\begin{array}
%[c]{c}%
%x_{1}^{\prime}\\
%x_{2}^{\prime}%
%\end{array}
%\right]  =P(t)\nabla f(x)=\left[
%\begin{array}
%[c]{c}%
%-2(x_{1}-1)(t+1)^{-2}\\
%-2(x_{2}-1)(t+1)^{-1}%
%\end{array}
%\right]  \text{, for }t\geq0\text{.}%
%\]
%
Then $\lambda_{1}(P(t))=(t+1)^{-2}$ and $\lambda_{2}
(P(t))=(t+1)^{-1}$, therefore $\int_{0}^{\infty}\lambda_{1}(P(t))dt<\infty$ while
$\int_{0}^{\infty}\lambda_{2}(P(t))dt=\infty.$ Thus, Theorem \ref{Th 1}.i
holds ($x_1(t)\equiv x_2(t) \equiv 1$ is a uniformly stable equibrium solution) while Theorem \ref{Th 1}.ii does not. The closed form solution
is given by $x_{1}(t)=1+c_{1}\exp(2/(t+1))$, $x_{2}(t)=1+c_{2}(t+1)^{-2}$ for
arbitrary constants $c_{1}$ and $c_{2}.$ We see that $\lim
_{t\rightarrow\infty}x_{1}(t)=1+c_1\neq1$; hence the equilibrium
solution ($x_1(t)\equiv1$) is not asymptotically stable.

More specifically, if $f$ has an isolated local maximum at $\overline{x}$ but
$\overline{x}$ is not an isolated critical point, then we can only conclude
stability, not asymptotic stability. 

\emph{Example 2.2 Asymptotic Stability requires Isolated Critical Point.} 
Here we produce an example which shows that the hypothesis H$_{2}$ is
essential to the conclusion that $\mathbf{x=}\overline{\mathbf{x}}$ is
asymptotically stable when $P(t)$ satisfies (EC) as in Theorem \ref{Th 1}.ii. Even for gradient systems, the necessity of adding the assumption that the point where the isolated local extremum occurs is also
an isolated critical point has been missed by some authors (e.g.,
Part 3 of the theorem on page 205 of Hirsch et al. \cite{hirsch_differential_2004}). 

We create a radially symmetric function $f(r,\theta)$ using polar coordinates that is
continuously differentiable on the unit circle $r\leq1$, that has an absolute
maximum value at the origin, that decreases as $r$ increases, and is such that
there is a sequence of concentric circles $r=r_{i}$ with $r_{i}$ decreasing to
zero as $i\rightarrow\infty$ and with each $r=r_{i}$ consisting entirely of
critical points of $f$. We again let $P(t)$ be the $2\times2$ identity matrix,
thus satisfying (EC). The system $\mathbf{x}^{\prime}=P(t)\nabla f$ will then
have the properties we seek, namely, we no longer have isolated critical
points of $f$.\newline\hspace*{0.25in}We first define sequences $x_{n}$ and
$z_{n}$ by $x_{n}=-2^{-n}$ and $z_{n}=(1-4^{-n})/3$ for $n=0,1,2,\cdots$. Let
$I_{n}$ be the interval $[x_{n},x_{n+1}]$. The union of the intervals $I_{n}$
is then the interval $[-1,0)$. We define a function $p(x)$ on the interval
$[-1,0)$ which restricted to the interval $I_{n}$ is a cubic polynomial
$p_{n}(x)$. Furthermore, we require each $p_{n}(x)$ to satisfy%
\begin{equation}
p_{n}(x_{n})=z_{n},\hspace*{0.1in}p_{n}^{\prime}(x_{n})=0,\hspace*{0.1in}%
p_{n}(x_{n+1})=z_{n+1},\hspace*{0.1in}p_{n}^{\prime}(x_{n+1})=0.\label{Req}%
\end{equation}

%\begin{figure*}
%% Use the relevant command to insert your figure file.
%% For example, with the graphicx package use
%  \includegraphics[width=1.0\textwidth]{Ex3.png}
%\caption{Graph of the bivariate polynomial spline used in Example 2.2 and the associated phase portrait. The dotted lines in the phase portrait indicate the locations the critical points given by concentric circles of radius $r_i=2^{-n}$.}
%% (only $n=0,1,2,3,4$ are shown)
%%The vector field clearly shows asymptotic stability will not occur because a trajectory starting outside $r=r_i$ is ``trapped'' away from $\bar{x}$ by the circle $r=r_i$ and cannot approach $\bar{x}$ as $t\rightarrow\infty$.
%%The plotted vectors clearly show asymptotic stability will not be achieved because the system will become ``trapped'' at next circle of smaller radius and the system is not guaranteed to approach $\bar{x}$ as $t\rightarrow\infty$. 
%%The smaller vectors indicate the diminished magnitude of change in the vicinity of critical points.
%\end{figure*}

\noindent Letting $p_{n}(x)=\alpha(x-x_{n})^{3}+\beta(x-x_{n})^{2}%
+\gamma(x-x_{n})+\delta$ and using the requirements in (\ref{Req}), we find
after some algebra and calculus that $\alpha=-2^{n+2}$, $\beta=3$, $\gamma=0$
and $\delta=(1-4^{-n})/3$. We then find, again using calculus, that
$p_{n}^{\prime}(x)>0$ for $x$ in the open interval $(x_{n},x_{n+1})$ and the
maximum value of $p_{n}^{\prime}$ on the interval $I_{n}$ is $3/2^{n+2}$. We
then extend $p(x)$ to the closed interval $[-1,0]$ by defining $p(0)=1/3$.
This makes $p$ continuous on $[-1,0]$. Considering difference quotients, it is
easy to see that the left-hand derivative of $p$ at $x=0$ exists and has value
zero. We symmetrically extend the definition of $p$ to the interval $[-1,1]$
by letting $p(x)=p(-x)$ for $0<x\leq1$. Taking into account the way the cubic
polynomials were pieced together at the endpoints and the fact that the
maximum value of $p^{\prime}(x)$ on the interval $I_{n}$ approaches zero as
$n\rightarrow\infty$, we see that $p$ has a continuous derivative on the
interval $[-1,1]$. 

Finally, we define the radially symmetric
$f(r,\theta)=f(r)$ by taking $f(r)=p(r)$ for $0\leq r\leq1$, $0\leq\theta
\leq2\pi$. Clearly, at any point on a circle $r=|x_{n}|$, we have
$f_{r}=f_{\theta}=0$ since $p^{\prime}(x_{n})=0$ and $f$ is independent of
$\theta$. Hence, all points on $r=|x_{n}|$ are critical points of $f$ and
yield equilibrium solutions of $\mathbf{x}^{\prime}=\mathbf{\nabla}f$. Even
though $f$ has an isolated maximum value at the origin, $x=(0,0)$ is not an
asymptotically stable equilibrium solution since solutions starting at $t=0$
between two concentric circles $r=|x_{n}|$ and $r=|x_{n+1}|$ are trapped in
that region and cannot approach the origin as $t\rightarrow\infty$. Of course,
Theorem \ref{Th 1}.i still applies to give that $x=(0,0)$ is a stable
equilibrium.  

\section{Basin of Attraction}

Given a uniformly asymptotically stable equilibrium $\overline{x}$ of (E), it
is of interest to know the set of points $x_{0}$ such that the trajectory
starting at point $x_{0}$ at some time $t_{0}$ exists for all $t\geq t_{0}$
and approaches $\overline{x}$ as $t$ tends to infinity; that is, the so-called
\textit{basin of attraction }of $\overline{x}.$ The following theorem provides
information on the basin of attraction in the setting of Theorem \ref{Th 1}.ii.

\begin{theorem}
[Basin of Attraction]\label{Th 2}Suppose H$_{0}$, H$_{1}$, H$_{2}$ and H$_{3}$
all hold. Let $M=f(\overline{x})$, let $c$ be a real number less than $M$ and
let $O_{c,\overline{x}}$ be the set defined by $O_{c,\overline{x}}%
=\{\overline{x}\}\cup\{x:c<f(x)<M\}$. Then $O_{c,\overline{x}}$ is open and
has a unique component $E_{c,\overline{x}}$ that contains $\overline{x}.$ Let
$\partial E_{c,\overline{x}}$ denote the boundary of $E_{c,\overline{x}}$, and
let $\overline{E_{c,\overline{x}}}=E_{c,\overline{x}}\cup\partial
E_{c,\overline{x}}$ denote the closure of $E_{c,\overline{x}}$. Consider
additional hypotheses: \hspace{0.1in}\newline\hspace*{0.3in}H$_{4}$:
$~E_{c,\overline{x}}$ is bounded and $\overline{E_{c,\overline{x}}}$ is
contained in $D$;\hspace{0.1in}\newline\hspace*{0.3in}H$_{5}$: $~f(x)=c$ for
all $x$ in $\partial E_{c,\overline{x}}$; and\hspace{0.1in}\newline%
\hspace*{0.3in}H$_{6}$:~ $f$ has no critical points other than $\overline{x}$
in $\overline{E_{c,\overline{x}}}$. \newline If H$_{4}$, H$_{5}$, and H$_{6}$
also hold, then $E_{c,\overline{x}}$ is contained in the basin of attraction
of $\overline{x}$.

\begin{proof}
Suppose \textit{H}$_{0}$ through \textit{H}$_{6}$ all hold. Let $M$, $c$ and
$O_{c,\overline{x}}$ be as defined above. We first prove $O_{c,\overline{x}}$
is open. By the continuity of $f$, the set $\{x:c<f(x)<M\}$ is open. Using the
continuity of $f$ and the fact that $f$ has an isolated maximum value at
$\overline{x}$, choose $\delta>0$ such that $f$ is defined on the ball
$B_{\delta}(\overline{x})$, $f(x)<M$ for $0<|x-\overline{x}|<\delta$ and
$|f(x)-M|<M-c$ for $x\in B_{\delta}(\overline{x})$. Then $c<f(x)<M$ for
$0<|x-\overline{x}|\,<\delta$ so $B_{\delta}(\overline{x})$ is open, contains
$\overline{x}$, and is contained in $O_{c,\overline{x}}$. It follows that the
set $O_{c,\overline{x}}=B_{\delta}(\overline{x})\cup\{x:c<f(x)<M\}$ is open.
Hence, there is a unique open component $E_{c,\overline{x}}$ of
$O_{c,\overline{x}}$ that contains the point $\overline{x}$.\newline%
\hspace*{0.25in}Let $x_{0}$ be any point of $E_{c,\overline{x}}$ and let
$x(t)$ be the solution of (E) satisfying the initial condition $x(t_{0}%
)=x_{0}$ for some $t_{0}\geq0$. We wish to prove that $x(t)$ exists for $t\geq
t_{0}$ and $\lim_{t\rightarrow\infty}x(t)=\overline{x}.$ This is clearly true
if $x_{0}=\overline{x}$ so we assume $x_{0}\neq\overline{x}$ and, in light of
the uniqueness of solutions to initial value problems, that $x(t)\neq
\overline{x}$ for all $t\geq t_{0}$. As before, we let $V(x)=M-f(x)$ for $x$
in $D$. While the trajectory $x(t)$ remains in $E_{c,\overline{x}}$, we have
by \textit{H}$_{6}$ that the trajectory derivative satisfies $V_{x}^{\prime
}(t)=-P(t)\nabla f(x(t))\cdot\nabla f(x(t))<0$. Because $V_{x}(t_{0})<M-c$ and
$V_{x}(t)$ decreases as $t$ increases, the trajectory $x(t)$ can never reach
$\partial E_{c,\overline{x}}$ where, by \textit{H}$_{5}$, $V_{x}(t)$ would
equal $M-c$. Hence, $x(t)$ stays in the region $E_{c,\overline{x}}$ and
therefore in the set $\overline{E_{c,\overline{x}}}$ so long as the solution
$x(t)$ continues to exist. By \textit{H}$_{4}$, $\overline{E_{c,\overline{x}}%
}$ is both closed and bounded and therefore compact. Since $x(t)$ stays in a
compact subset of $D$, it follows directly from Theorem 3.1 of Hartman
\cite{hartman_ordinary_1964}, that the right-maximal interval of existence of $x(t)$
as a solution of (E) cannot be of the form $[t_{0},\omega)$ with
$\omega<\infty$. Thus, the solution $x(t)$ exists for all $t\geq t_{0}%
$.\newline\hspace*{0.25in}From here on, the proof essentially follows that of
Theorem \ref{Th 1}.ii, but we repeat some of the details for clarity. First,
let $\lim_{t\rightarrow\infty}V_{x}(t)=\alpha$ and suppose $\alpha>0$. Then
using the continuity of $V$, find $\delta>0$ small enough that $0<V(x)<\alpha$
for $0<|x-\overline{x}|<\delta$. Now the set $\overline{E_{c,\overline{x}}%
}\backslash B_{\delta}(\overline{x})$ is closed and bounded by \textit{H}%
$_{4}$, so, by \textit{H}$_{6}$, the continuous function $\nabla
f(x)\cdot\nabla f(x)$ assumes a positive minimum $m_{1}$ on that set. Because
$x(t)$ never enters the set $B_{\delta}(\overline{x})$ where we would have
$V_{x}(t)=V(x(t))<\alpha$, we get that $V_{x}^{\prime}(t)\leq-m_{1}\lambda
_{1}(P(t))$ for $t\geq t_{0}$. This leads to $V_{x}(t)\rightarrow-\infty$ as
$t\rightarrow\infty$, a contradiction which shows that $\alpha=0$%
.\newline\hspace*{0.25in}The next step is to prove that $\lim_{t\rightarrow
\infty}x(t)=\overline{x}$. To do this, suppose $\lim_{t\rightarrow\infty
}x(t)\neq\overline{x}$. There then exists an $\varepsilon>0$ such that
$|x(t)-\overline{x}|\geq\varepsilon$ for arbitrarily large values of $t$. The
function $V(x)$ is positive and continuous on the compact set $\overline
{E_{c,\overline{x}}}\backslash B_{\varepsilon}(\overline{x})$, hence, $V(x)$
has a positive mimimum, call it $m_{2}$, on the set $\overline{E_{c,\overline
{x}}}\backslash B_{\varepsilon}(\overline{x})$. However, there are arbitrarily
large values of $t$ where $x(t)\in\overline{E_{c,\overline{x}}}\backslash
B_{\varepsilon}(\overline{x})$ for which $V_{x}(t)=V(x(t))\geq m_{2}$. This
contradicts $\lim_{t\rightarrow\infty}V_{x}(t)=0$ and completes the proof.
\end{proof}
\end{theorem}

We note that LaSalle's Theorem can be used to obtain information on the basin
of attraction of an equilibrium solution---for examples see Theorem 6.1 in
Leighton \cite{leighton_introduction_1976}, Theorem 11.11 in Miller and Michel
\cite{miller_ordinary_1982}, or the theorem on p. 200 of Hirsch et al.
\cite{hirsch_differential_2004}. However, those results deal with autonomous systems and do
not apply to (E).

\emph{Example 3.1 Basins of Attraction.} We conclude by giving an example illustrating both the use of Theorem
\ref{Th 2} and the role played by hypotheses \textit{H}$_{5}$ and
\textit{H}$_{6}$ of that theorem. Let $f(x_{1},x_{2})=96x_{2}-84x_{2}^{2}+28x_{2}%
^{3}-3x_{2}^{4}-10(x_{1}-2)^{2}$. Let $P(t)$ be any $2\times2$ matrix-valued
function defined and continuous for $t\geq0$ and such that the eigenvalue
condition (EC) holds. Equation (E) becomes%

%\begin{figure}[!b]
%% Use the relevant command to insert your figure file.
%% For example, with the graphicx package use
%  \includegraphics[width=1\textwidth]{Example7A.png}
%\caption{Plot of $f(x_1,x_2)$ from Example 3.1. In the plot on the left, the level sets $L_{37}$ (black), $L_{33}$ (red), and $L_{20}$ (blue) are superimposed on the surface. The plot on the right depicts the sets of points interior to the level set curves using the same color scheme (e.g., $E_{20,p_1}$  is blue like $L_{20}$).} 
%%For illustration purposes, the elevation of the set was chosen to be the number $c$ in 
%%the definition of  $E_{c,p_i}$; thus the labels for each set include the cross-product with the elevation (e.g., $E_{20,p_1}\times\{20\}$). Theorem \ref{Th 2} does not apply to  $E_{20,p_1}$ because the function $f$ assumes the value 37 on the inside boundary and the value 20 on the outside boundary, thus violating hypothesis $H_5$.
%\end{figure}

\[
\left[
\begin{array}
[c]{c}%
x_{1}^{\prime}\\
x_{2}^{\prime}%
\end{array}
\right]  =P(t)\left[
\begin{array}
[c]{c}%
-20(x_{1}-2)\\
-12(x_{2}-1)(x_{2}-2)(x_{2}-4)
\end{array}
\right]  \text{.}%
\]

\noindent Then $f$ has local maximum values at the points $p_{1}=(2,1)$ and
$p_{2}=(2,4)$ and a saddle at $p_{3}=(2,2)$ with $f(2,1)=37$, $f(2,4)=64$ and
$f(2,2)=32$. 

For a real
number $c$, let $L_{c}$ denote the level set defined by $L_{c}=\{(x_{1}%
,x_{2}):f(x_{1},x_{2})=c\}$. 
$L_{33}$ consists of two simple closed curves; we let $C_{p_{1}}$ and
$C_{p_{2}}$ denote the curve having the point $p_{1}$ and $p_{2}$
(respectively) as an interior point. Then the set $E_{33,p_{1}}$consists of
all points interior to $C_{p_{1}}$ while $E_{33,p_{2}}$ consists of all points
interior to $C_{p_{2}}$. Theorem \ref{Th 2} applies and shows that all
trajectories $x(t)$ having $x(t_{0})$ in $E_{33,p_{1}}$ tend to $p_{1}$ as
$t\rightarrow\infty$, with a similar conclusion for trajectories in
$E_{33,p_{2}}$. It is interesting to consider $E_{20,p_{1}}=\{p_{1}%
\}\cup\{x:20<f(x)<37\}$ and $E_{20,p_{2}}=\{p_{2}\}\cup\{x:20<f(x)<64\}$.
First, $E_{20,p_{2}}$ contains all points interior to a simple closed curve
containing both $p_{1}$ and $p_{2}$ in its interior; hence, Theorem \ref{Th 2}
does not apply to $E_{20,p_{2}}$ because \textit{H}$_{6}$ is violated. On the
other hand, $E_{20,p_{1}}$ consists of $E_{20,p_{2}}\backslash\overline
{E_{37,p_{2}}}$. Now, Theorem \ref{Th 2} does not apply to $E_{20,p_{1}}$
because the boundary of $E_{20,p_{1}}$ contains points of the level set
$L_{37}$ at which $f$ takes on the value $37$ thus violating \textit{H}$_{5}$;
clearly some trajectories starting in $E_{20,p_{1}}$ will tend toward the
boundary points in $L_{37}$ while others will tend toward $p_{1}$.

%\section*{Conclusions}
%
%We have demonstrated conditions for the existence and stability of equilbria
%in modified gradient systems. Such systems include the continuous-time multivariate breeders
%equation with a time-dependent $\mathbf{G}$-matrix which is frequently used by evolutionary biologists. Our research shows that biologists can simply search for the
%isolated local maxima of a fitness function; these points are guaranteed to be
%at least uniformly stable. Furthermore, in cases where the smallest eigenvalue
%of the G-matrix, $\lambda_{1}(P(t))$ in our notation, meets the condition $\int_{0}^{\infty
%}\lambda_{1}(P(t))dt=\infty$, then the equilibrium is guaranteed to be
%uniformly asymptotically stable. Finally, an understanding of the inverse
%image under the fitness function $f$ of intervals of the form $(c,f(\overline
%{x}))$ gives information on the basin of attraction of an equilibrium at
%$\overline{x}$.

%I included the below from another AML article to check the formatting of the documenting so that we
%know how long our article is with regard to AML 6pg limit.
%Notice, the second order derivative of σn arises in (3.17), which means our algorithm uses a smooth function to approximate the exact solution in admissible set $\Sigma$. The numerical implementations of this algorithm will be given elsewhere.

\section*{References}
 
\bibliography{EvolutionaryStability}

\end{document}